\documentclass[reqno]{amsart}
\usepackage{amscd}
\usepackage{amssymb}
\usepackage{verbatim}
\usepackage{emlines2}
\usepackage{bezier}
\usepackage{epsfig}

\usepackage{fancyheadings}
\theoremstyle{plain}
\newtheorem{Thm}{Theorem}[section]

\newtheorem{Def}[Thm]{Definition}
\newtheorem{Lem}[Thm]{Lemma}

\newtheorem{Rem}[Thm]{Remark}
\newtheorem{Main}[Thm]{Main Theorem}

\newtheorem{LBT}[Thm]{Lower Bound Theorem}
\newtheorem{UBT}[Thm]{Upper Bound Theorem}

\begin{document}
\parindent0mm
\newcommand{\Aut}{\ensuremath{\text{\normalfont Aut}}}
\renewcommand{\labelenumi}{(\roman{enumi})} 
\title[A classification of $12$-vertex triangulations of $S^2 \times S^2$]{A 
classification of centrally-symmetric and cyclic $12$-vertex triangulations of 
$S^2 \times S^2$}
\thanks{The second author was supported by the Deutsche Forschungsgemeinschaft
(DFG, Grant \linebreak Ku 1203/2-1)}
\maketitle
Gunter Lassmann$^1$ and Eric Sparla$^2$\\


\small

$^1$Forschungs- und Technologiezentrum der Deutschen Telekom (TZD FE 124),\\
Goslarer Ufer 35, 10589 Berlin, Germany\\
lassmann@tzd.telekom.de\\[2mm]

$^2$Mathematisches Institut B4, Universit\"at Stuttgart,\\
Pfaffenwaldring 57, 70569 Stuttgart, Germany\\
sparla@mathematik.uni-stuttgart.de
\normalsize

\begin{abstract} 
In this paper our main result states that there exist exactly three 
combinatorially distinct centrally-symmetric $12$-vertex-triangulations 
of the product of
two $2$-spheres with a cyclic symmetry. We also compute the automorphism 
groups of the triangulations. These instances suggest that there
is a triangulation of $S^2 \times S^2$  with $11$ vertices -- the 
minimum number of vertices required.\\ 
\end{abstract}


\section{Introduction and Main Theorem} \label{S:intro}
By a theorem of K\"uhnel \cite{Ku3} any triangulation
of $S^2 \times S^2$ must have at least $10$ vertices and any triangulation 
with only $10$ vertices must be $3$-neighborly, i.~e.\
any $3$-tuple of vertices spans a $2$-simplex of the triangulation.
By a computer-aided enumeration K\"uhnel and the first author \cite{K-L}
proved that there does not exist a $3$-neighborly triangulation of a
$4$-manifold with $10$ vertices, hence any 
triangulation of $S^2 \times S^2$ requires $11$ vertices or more. 
Up to now no triangulation of $S^2 \times S^2$ with only $11$ vertices has been
published (cf.\ however Remark $4.2$), but the second author has found in 
\cite{Sp1} a highly-symmetric
$12$-vertex-triangulation of $S^2 \times S^2$. In this paper we present a
classification of all centrally-symmetric $12$-vertex-triangulations of 
$S^2 \times S^2$ with a cyclic symmetry. 
We will establish that one of the triangulations of this
classification is combinatorially equivalent to the one found in \cite{Sp1}.

In order to state this main result we need the notion of combinatorial 
manifolds as follows:
 
\begin{Def}
A simplicial complex $M$ is called a {\em combinatorial 
manifold of dimension n} if the dimension of any simplex of $M$ is at
most $n$ and if  the link $lk(\Delta^k,M) :=lk(\Delta^k)$ of any $k$-simplex
$\Delta^k$ of $M$, $0 \leqslant k \leqslant n$, is a triangulated
$(n-k-1)$-sphere.

A diagonal is an edge consisting of vertices of $M$ that itself does 
not belong to $M$.
 
Two combinatorial manifolds (or more generally two simplicial complexes) $M_1$ 
and $M_2$ are said to be combinatorially equivalent, in this case we write 
$M_1 \sim M_2$,
if there is a bijective mapping $\phi$ between the faces of $M_1$ and $M_2$ 
that is inclusion-preserving, i.~e.\ such that $F_1 \subset F_2$ if and only 
if $\phi(F_1) \subset \phi(F_2)$.

The underlying complex of a combinatorial manifold (or a simplicial complex)
$M$, written $|M|$, is the union of all simplices of $M$.
\end{Def}

In the sequel the vertices of all examples will be denoted by
$\langle i \rangle$ with $0 \leqslant i \leqslant 9$, $\langle a \rangle$
for the 11th and $\langle b \rangle$ for the 12th vertex. 
If the permutation $\zeta := (0, 1, \ldots, 9, a, b)$ is an automorphism of
the triangulation we regard, we talk about $\zeta$ as a cyclic symmetry 
generating the cyclic group $C_{12}$. Using the usual notation $f_i(M)$ 
(or simply $f_i$ if there is no danger of confusion) for the number of 
$i$-faces of a combinatorial manifold $M$, $\chi(M)$ for its 
Euler-characteristic and the above definition we are now able to state our 

\begin{Main}
There exist exactly three combinatorially distinct types of combinatorial
$4$-manifolds $M_i$, $1 \leqslant i \leqslant 3$, with $\zeta$ as an
automorphism such that the $1$-simplex $\langle 0 6 \rangle$ is a 
diagonal of $M_i$ and such that $f_4(M_i)=72=6f_0(M_i)=18\chi(M_i)$. 
For all the $M_i$ it holds $|M_i| \simeq_{PL} S^2 \times S^2$, i.~e.\ the
$M_i$ are $PL$-homeomorphic to $S^2 \times S^2$, and their automorphism
groups $\Aut(M_i)$ are as follows:
\[ \Aut(M_i) \; \cong \; 
\begin{cases} 
    C_{12} \rtimes C_2 & \text{ for }\: i=1 \\
    A_5 \rtimes C_4  & \text{ for }\:i=2 \\
    C_{12} \rtimes (C_2 \times C_2) & 
    \text{ for }\: i=3
\end{cases}
\]
where $\rtimes$ denotes the semidirect product and $C_k$ the cyclic group
of order $k$.\\
\end{Main}

Before proving this result in Section $2$, we will now describe several 
interesting properties that are common to all $M_i$:

The $f$-vector of $M_i$ is uniquely determined by the 
Dehn-Sommerville equations \cite{Ku2} for triangulated $4$-manifolds
\[ 2f_1-3f_2+4f_3-5f_4 = 0 = 2f_3 - 5f_4 \]
to be
\[ (f_0,f_1,f_2,f_3,f_4)=(12,60,160,180,72).\]
\vspace{0.5mm}\relax{}\\
The cyclic symmetry immediately implies that the automorphism group of $M_i$
is transitive on its vertices. This proves that the
vertex links of all vertices are combinatorially equivalent.\\

The cyclic symmetry also yields in connection with the condition
that $\langle 0 6 \rangle$ is a diagonal of $M_i$, 
that $M_i$ has at least $6$ pairwise disjoint diagonals. Hence
$M_i$ can be embedded in the $6$-dimensional cross-polytope $C_6^\ast$.
As $f_i(M_i)=f_i(C_6^\ast)=2^{i+1} \binom{6}{i}$ for $i=0,1,2$ any such
embedding is $2$-Hamiltonian, i.~e.\ contains the $2$-skeleton of 
$C_6^\ast$, and therefore $M_i$ is simply-connected \cite[$3.8$]{Ku1}.
\pagebreak
From these properties we conclude that all $M_i$ satisfy the assumptions 
of the following Lower Bound Theorem. Moreover they are examples for the
case of equality in the inequality, hence they prove that the inequality
is sharp:

\begin{LBT}[{\cite[$1.2$]{Sp2}}]  Assume $M$ in $\mathbb{E}^d$ is a
combinatorial $2k$-manifold, whose convex hull $\mathcal{H}(M) = P$ is
-- up to affine transformations -- a centrally-symmetric simplicial 
$d$-polytope $P \subset \mathbb{E}^d$. Let $M$ contain the $k$-skeleton
$Sk_k(P)$ of $P$, that is the set of all faces of $P$ of dimension at most 
$k$, and let $M$ be a subcomplex of the boundary complex 
$\mathcal{C}(\partial P)$ of $P$. Then the following statements hold:
\begin{enumerate}
\item $(-1)^k\binom{2k+1}{k+1} (\chi(M)-2) \geqslant 4^{k+1}
\binom{\frac{1}{2}(d-1)}{k+1}$.  
\item For $d>2k+1$ equality in $(1)$ is attained if and only if $P$ is 
affinely equivalent to the $d$-dimensional cross-polytope $C_d^\ast$.
\end{enumerate} 
\vspace{2mm}
\end{LBT}

The $M_i$ also satisfy equality in the next Upper Bound Theorem. A required
fixed-point free involution is clearly given by $\zeta^6$.

\begin{UBT}[{\cite[4.1]{Sp1}}]  Let $M$ be a combinatorial $4$-manifold
 with $n$ vertices and with a fixed-point involution $\sigma$ acting on the 
 triangulation. Then $n=2m$ is even and the inequality
 \[ 10(\chi(M)-2) \leqslant \frac{4}{3}(m-1)(m-3)(m-5) = 4^3\binom
{\frac{1}{2}(m-1)}{3}.\]
 holds with equality if and only if $M$ can be embedded in $C^\ast_m$ such
 that this embedding contains the $2$-skeleton of $C^\ast_m$.\\ 
 \end{UBT}


For the proof of the Main Theorem we are going to describe explicitly in the 
next section the triangulations $M_i$. From the lists given there it becomes 
obvious that they are all centrally-symmetric.
In Section $3$ we state further properties of the $M_i$. 
Finally in Section $4$ we conclude with another interesting 
$12$-vertex triangulation found while  proving 
our result and we make some remarks concerned with $11$-vertex triangulations 
of $S^2 \times S^2$.\\[2mm]

\section{Proof of the Main Theorem}
\texttt{SUNI}, a computer program written by the first author and described 
in detail in \cite{L1}, is used for 
determining all possible candidates for $4$-manifolds that satisfy 
several conditions. For instance all candidates found by \texttt{SUNI} are 
cyclic and satisfy the Dehn-Sommerville equations. Additionally a 
combinatorial test checks if the
Euler-charateristic of all edge-links equals $2$. In the case of our 
additional assumptions -- in particular $\langle 0 6 \rangle$ has to be a 
diagonal -- the program delivers three possible candidates:\\

\renewcommand{\labelenumi}{(\arabic{enumi})} 
\begin{enumerate}
\item
$ \quad 1\, 1\, 1\, 1\,8 \qquad 2\,2\,3\,2\,3 \qquad 1\,1\,2\,1\,7 \qquad
1\,4\,3\,2\,2 \qquad 1\,1\,3\,5\,2 \qquad 1\,3\,1\,3\,4$ 
\item 
$ \quad 1\, 1\, 1\, 1\,8 \qquad 2\,2\,3\,2\,3 \qquad 1\,1\,2\,1\,7 \qquad
1\,4\,3\,2\,2 \qquad 1\,1\,3\,5\,2 \qquad 1\,3\,1\,4\,3$ 
\item 
$ \quad 1\, 1\, 1\, 1\,8 \qquad 2\,2\,3\,2\,3 \qquad 1\,1\,2\,5\,3 \qquad
1\,1\,3\,5\,2 \qquad 1\,1\,3\,4\,3 \qquad 1\,2\,1\,4\,4$
\vspace{2mm}
\end{enumerate}
\renewcommand{\labelenumi}{(\roman{enumi})} 

The notation $(y_1,\ldots,y_d)$ (or shorter $\;y_1 \,\ldots\,y_d\;$) denotes
a difference $d$-cycle, that generates a $C_n$-orbit of $d$-simplices
as follows:
\[ \left\{ \langle x \; x+y_1 \; \ldots \; \left. 
x+\sum\limits_{i=1}^{d-1} y_i \rangle 
\; \right| \; x \in C_n \right\}. \]
(Note that  $y_d := n-\sum_{i=1}^{d-1}y_i$.)
In the case of Theorem $1.2$ we have $d=5$ and $n=f_0=12$.
Let us denote candidate $(i)$ by $M_i$ and for the sake of brevity let 
us omit the symbols ``$\langle$'' and ``$\rangle$'' in the list of its
simplices, then the complete list 
of $4$-simplices of, for example, $M_1$ is given by

\[ \begin{array}{lllllllllll}
01234\; && 02479\; && 01245\; && 0158a\; && 0125a\; && 01458 \\[1mm]
12345\; && 1358a\; && 12356\; && 1269b\; && 1236b\; && 12569 \\[1mm]
23456\; && 2469b\; && 23467\; && 237a0\; && 23470\; && 2367a \\[1mm]
34567\; && 357a0\; && 34578\; && 348b1\; && 34581\; && 3478b \\[1mm]
45678\; && 468b1\; && 45689\; && 45902\; && 45692\; && 45890 \\[1mm]
56789\; && 57902\; && 5679a\; && 56a13\; && 567a3\; && 569a1 \\[1mm]
6789a\; && 68a13\; && 678ab\; && 67b24\; && 678b4\; && 67ab2 \\[1mm]
789ab\; && 79b24\; && 789b0\; && 78035\; && 78905\; && 78b03 \\[1mm]
89ab0\; && 8a035\; && 89a01\; && 89146\; && 89a16\; && 89014 \\[1mm]
9ab01\; && 9b146\; && 9ab12\; && 9a257\; && 9ab27\; && 9a125 \\[1mm] 
ab012\; && a0257\; && ab023\; && ab368\; && ab038\; && ab236 \\[1mm]
b0123\; && b1368\; && b0134\; && b0479\; && b0149\; && b0347
\vspace{2mm}
\end{array} \]

The simplices of the other candidates can be calculated analogously from 
the difference-$5$-cycles. Obviously the $M_i$ contain all edges except
the diagonals \linebreak $\langle j (j+6) \rangle$, where $j$ as well as 
$j+6$ are regarded $\mod 12$. This is consistent with the theoretical 
value $f_1(M_i)=60$.

All candidates are centrally-symmetric as $\langle a_1,\ldots, a_j
\rangle \in M_i$ implies $\langle b_1,\ldots,b_j\rangle \in M_i$ where 
$b_l = a_l +6 \mod 12$, $1\leqslant l \leqslant j$.\\

It turns out that $M_2$ is combinatorially equivalent to the
example $M$ found in \cite[Theorem $3.1$]{Sp1}. The images of the vertices
of $M$ (in the notation used in \cite{Sp1}) of a required bijection
$\phi : M \mapsto M_2$ are as follows: \label{M_old}
\[ \begin{array}{cccccc} \phi(0)=5, & \phi(1)=2, & \phi(2)=0, &
 \phi(3)=7, &  \phi(4)=3, &  \phi(5)=4, \\
\phi(\bar{0})=b, &\phi(\bar{1})=8, & \phi(\bar{2})=6, & \phi(\bar{3})=1, &
\phi(\bar{4})=9, & \phi(\bar{5})=a.  
\end{array}\]
Unfortunately, the group $A_5 \times C_2$ discussed in \cite{Sp1} is not the 
full automorphism group of $M$ but only a subgroup of index $2$ of it.  
The group given here is its full automorphism group.\\ 

For proving the Main Theorem it suffices to show the following Lemmata. 
The first two lemmata establish the classification and are therefore the most 
important part of the proof of $1.2$. In the third lemma we calculate
the automorphism groups.


\begin{Lem}
$M_1 \not\sim M_2 \not\sim M_3 \not\sim M_1$, i.~e.\ there
are three mutually combinatorial inequivalent triangulations.
\end{Lem}

\begin{Lem}
$M_i$, $1 \leqslant i \leqslant 3$, are combinatorial
$4$-manifolds with $|M_i| \simeq_{PL} S^2 \times S^2$.
\end{Lem}

\begin{Lem}
$\Aut(M_i) \; \cong \; 
\begin{cases} 
    C_{12} \rtimes C_2 & \text{ for }\: i=1 \\
    A_5 \rtimes C_4  & \text{ for }\:i=2 \\
    C_{12} \rtimes (C_2 \times C_2) & 
    \text{ for }\: i=3
\end{cases}\;.$
\vspace{5mm}
\end{Lem}

\begin{proof}[Proof of Lemma $2.1$]
The links of all edges of $M_2$ are subdivided octahedra as
sketched in figure $1$.
For the other $M_i$ this is not the case. Therefore $M_2 \not\sim M_v$
for $v\in \{1,3\}$.\\
\begin{center}
\epsfig{file=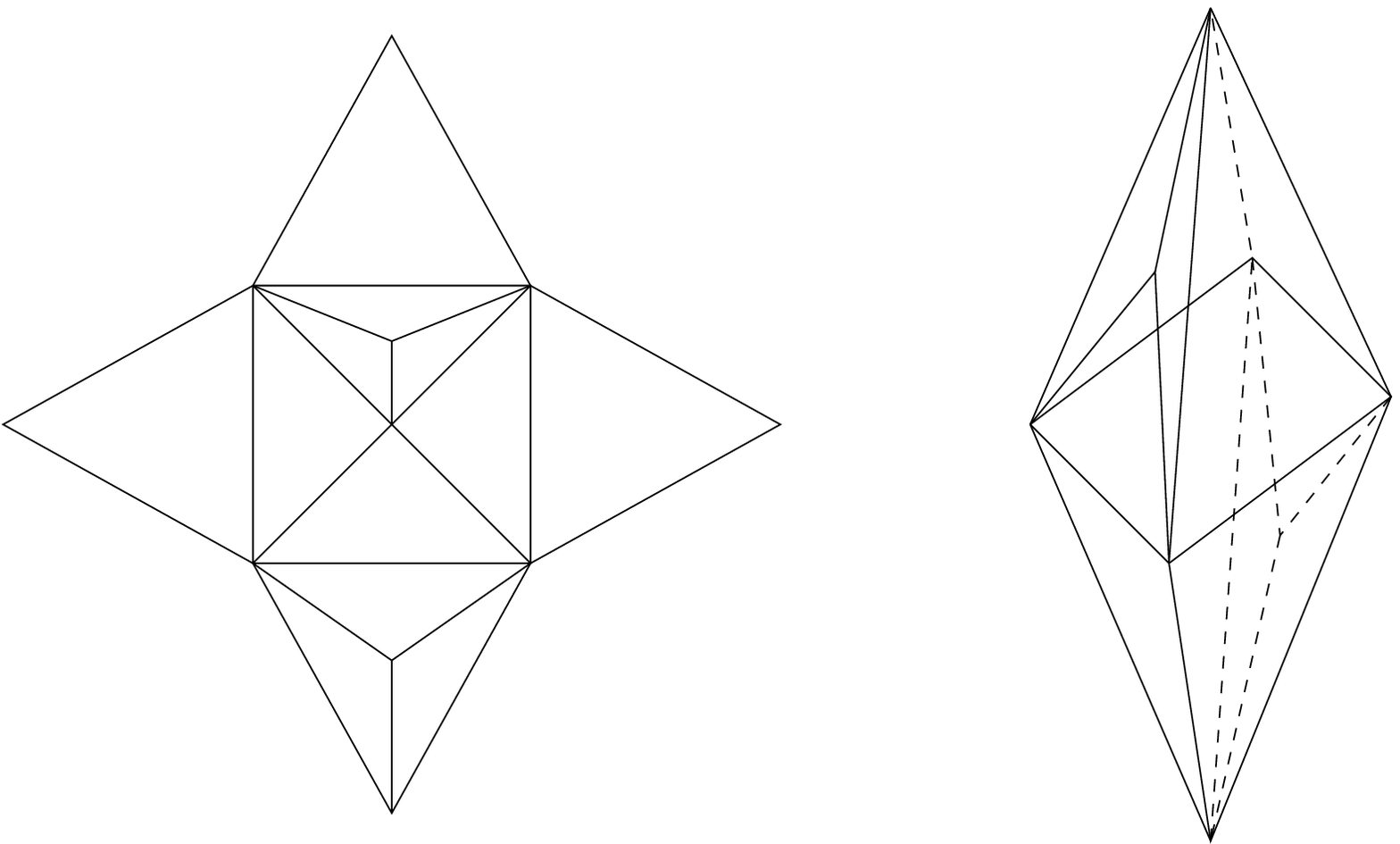, width=9.4cm, height=5.2cm} \\
\end{center}
\begin{tabular}{rl} 
{\bf Figure 1.} & On the left-hand side the link of an edge of $M_2$, $M_4$
respectively is \\ & sketched and on the right-hand side 
its $3$-dimensional 
realisation. 
\end{tabular}\\[4mm]
In the link $lk(k,M_i)$ of each edge $k$ of $M_i$ every vertex is 
adjacent to $j$ other vertices where $3 \leqslant j \leqslant 6$. 
The number of vertices with $j$ adjacent vertices in $lk(k,M_i)$ is
denoted by $val(j)$ and
\[ v(lk(k,M_i)) := (val(3),val(4),val(5),val(6)) \]
is called the valence-vector of $lk(k,M_i)$. As the links of all vertices are 
combinatorially equivalent, it can only hold $M_{i_1} \sim
M_{i_2}$, if for all $1 \leqslant j_1 \leqslant 5$ there exists a 
$j_2$ with $1 \leqslant j_1 \leqslant 5$ such that 
\[ v(\langle 0j_1 \rangle, M_{i_1}) = v(\langle 0j_2 \rangle, M_{i_2}).\]
But we get 
\[ v(\langle 03 \rangle, M_3) = (0,4,4,0) \neq 
v(\langle 0j \rangle, M_1) = \begin{cases}
 (1,3,3,1) & \text{ for } j=1 \\
 (2,0,6,0) & \text{ for } j=2 \\
 (2,2,2,2) & \text{ for } j=3\\
 (1,4,1,2) & \text{ for } j=4 \\
 (1,3,3,1) & \text{ for } j=5 
\end{cases}
\]
Hence $M_1 \not\sim M_3$.
\end{proof}
\pagebreak

\begin{proof}[Proof of Lemma $2.2$]
As we have shown the statement for $M_2$ already in \cite{Sp1} 
(up to combinatorial equivalence, cf.\ page \pageref{M_old}) and as the proofs for
$M_1$ and $M_3$ are similar, but not completely analogous to the one given in
\cite{Sp1} we will present here a proof only for $M_3$.\\

As $\Aut(M_3)$ operates transitively on the vertices if suffices to show
that the link $lk(0)$ of vertex $0$ is a triangulated $S^3$. 
To do this we have to establish:
\begin{enumerate}
\item 
The link of each vertex of $lk(0)$, i.~e.\ the link of each edge containing
vertex $\langle 0 \rangle$, is a triangulated $2$-sphere. Because
of automorphism $\zeta$ we only have to check this for the vertices $\langle
i \rangle$ with $1 \leqslant i \leqslant 5$.
Here we can sketch the triangulations in the same way as shown in Figure $1$
and we get $2$-spheres with the valence vector we have calculated in the
preceding lemma.
\item
$lk(0)$ is a triangulated $3$-sphere. To show this we split up the vertex 
set of $lk(0)$ into two disjoint subset and regard their span:
\[ A:=span(1,2,3,4,7) \text{ and } B:=span(5,8,9,a,b).\] 
$A$ and $B$ each consist of two tetrahedra and are therefore collapsible. 
Hence $A$ and $B$ are $3$-balls (cf.\ \cite{R-S}).
Cluing these balls together proves that the underlying complex of $lk(0)$ is
homeomorphic to $S^3$.\\ 
\end{enumerate}

It remains to determine the topological type of $|M_3|$. Let
\[ \alpha_1 := span(0,2,4,6,8,a)\;\text{ and }\;\alpha_2 := span(0,1,3,5). \]
Then $\alpha_1 \sim \partial C_3^\ast$ and $\alpha_2 \sim \Delta^3$. 
It turns out that $\alpha_1$ together with $\alpha_2$ generate the second 
homology group $H_2(M_3,\mathbb{Z})$ and that the intersection form of $M_3$ 
equals the intersection form of $S^2 \times S^2$, namely
\[ \begin{pmatrix} \alpha_1\cdot \alpha_1 & \alpha_1 \cdot \alpha_2 \\
\alpha_2\cdot \alpha_1 & \alpha_2 \cdot \alpha_2 \end{pmatrix} = 
\begin{pmatrix} 0 & 1 \\ 1 & 0 \end{pmatrix}.\]
This is an even form and by Freedman's Theorem \cite{Fr} there exists
exactly one simply-connected, closed, topological $4$-manifold 
representing that form. Hence $|M_3|$ must be homeomorphic 
to $S^2 \times S^2$.\\ 
(To get for instance the intersection number $\alpha_1 \cdot \alpha_2$ 
regard 
\[ N_1 := span(0,1,3,5) \text{ and } N_2 := span(0,2,4,6,7,8,9,a,b).\]
As
\[ lk(0) \setminus \text{span}(0,1,3,5) = lk(0) \setminus 
(N_1 \cap lk(0)) \searrow N_2 \cap lk(0)\]
(where ``$X \searrow Y$'' means ``$X$ collapses onto $Y$'' \cite{R-S}),
the complement of $N_1 \cap lk(0)$ is homotopically equivalent to $N_2 \cap 
lk(0)$. Because
\[ |lk(0) \setminus \text{span}(0,1,3,5)|=\langle 13\rangle \cup \langle 35 
\rangle \cup \langle 51\rangle\] and
\[ |N_2 \cap lk(0)| \searrow \langle 24\rangle \cup \langle 48 \rangle
\cup \langle 8 a\rangle \cup \langle a2 \rangle\]
we get two disjoint cycles that must be unknotted and linked in $lk(0)$. 
Consequently $\alpha_1 \cdot \alpha_2 = \pm 1$ and with an appropriate
orientation we get the intersection number stated above.)

To prove that $|M_3|$ is even $PL$-homeomorphic to $S^2 \times S^2$ we can
mimic the prove given in \cite{Sp1} that essentially uses that there is only
one $3$-ball bundle over $S^2$ with the intersection form of $S^2 \times S^2$.
This proof does not use Freedman's Theorem but we still need to calculate
the intersection form.
\end{proof}
\vspace{2mm}

\begin{proof}[Proof of Lemma $2.3$]
We use a simple \texttt{GAP}-program\footnote{More information on the 
freeware-algebra-system GAP can be found under the URL 
\texttt{http://www-gap.dcs.st-and.ac.uk/\symbol{126}gap/}.
For information on the GAP-program used here, send an E-Mail to the second
author.} written by the second author. This program explicitly
calculates all bijections between the vertex-links of two
candidates. The order of $\Aut(lk(0,M_i))$ is 
$2$ for $i=1$ and $4$ for $i=3$. The automorphisms of these cases are
\[ id,\;\alpha=(1,5)(2,a)(4,8)(7,b) \quad 
\text{for} \quad i=1\]
and 
\[ id, \;  \beta_1=(1,b)(2,a)(3,9)(4,8)(5,7), \; \beta_2=(1,7)(3,9)(5,b),\]
\[ \beta_3 = (1,5)(2,a)(4,8)(7,b) \quad \text{for} \quad i=3. \]
For $i=2$ we get $20$ automorphisms. In all cases it is easy to show
\[ \Aut(lk(0, M_i)) = \begin{cases} C_2 & \text{for } i=1 \\
C_5 \rtimes C_4 & \text{for } i=2 \\
C_2 \times C_2 = A^{(2,2)} & \text{for } i=3\end{cases},\]
where $A^{(2,2)}$ denotes the Kleinian group.

(Note that $\Aut(lk(0,M_2))$ is transitive on the vertices of $lk(0,M_2)$ 
because
\[ (2,a,b,9,7)(1,8,4,5,3) \in \Aut(lk(0,M_2)) \ni (1,a,8,9)(2,3,7,4)(5,b).\]
By the cyclic symmetry $\Aut(M_2)$ acts therefore transitively on the edges of 
$M_2$. This is not the case for $i=1$ and $i=3$.)

The transitivity of $\Aut(M_i)$ on the vertices leads to 
\begin{equation} 
|\Aut(M_i)| \leqslant f_0(M_i)|\Aut(lk(0,M_i)| = 
\begin{cases} 
24 & \text{for } i=1\\ 
240 & \text{for } i=2\\
48 & \text{f\"ur } i=3
\end{cases}.
\end{equation}
Denote by $A \cdot B$ the complex product of two groups $A$ and $B$.\\
As (i) $|\langle \zeta \rangle \cdot \langle \alpha \rangle| \cong
|C_{12}\cdot C_2| = 24$, (ii) $\alpha\in \Aut(M_1)$ and
(iii) $C_{12} \triangleleft C_{12} \cdot C_2$, inequality $(1)$ implies that 
$C_{12} \cdot C_2 \cong C_{12} \rtimes C_2$ 
is the full automorphism group of $M_1$.\\

Analogously it can be shown $\Aut(M_3) \cong C_{12} \cdot \langle \beta_1,
\beta_2,\beta_3\rangle$. $\Aut(lk(0,M_3))$ is not normal in $\Aut(M_3)$ as 
\[ \zeta\cdot\langle \Aut(lk(0,M_3))\rangle\neq\langle \Aut(lk(0,M_3))
\rangle\cdot\zeta,\]
but 
\[ \beta_i \langle \zeta\rangle  = \langle \zeta\rangle\beta_i \quad
\text{for}\quad 1 \leqslant i \leqslant 3\]
and therefore $\langle \zeta \rangle \triangleleft \Aut(M_3)$.
Hence $\Aut(M_3) \cong C_{12} \rtimes A^{(2,2)}$.\\

It remains to discuss the automorphism group of $M_2$. $A_5 \times C_2$ is a
subgroup of $\Aut(M_2)$ as already seen in \cite{Sp1}.
This subgroup is generated by the permutations 
\[ (0,7,3,4,2)(1,9,a,8,6), \quad (0,6)(1,7)(2,8)(3,9)(4,a)(5,b)\vspace{-2mm}\]
\[\text{and }\quad\quad (0,6)(1,9)(2,5)(3,7)(4,a)(8,b). \]
(Here we had to rename the vertices according to the bijection $\phi$ stated
above.) Additionally $\zeta$ is a (cyclic) automorphism operating on $M_2$.
According to \texttt{GAP} the order of the group generated by these four
permutations is $240$ and therefore by $(1)$ it must be the full
automorphism group. A closer examination shows \linebreak$\Aut(M_2) \cong
A_5 \rtimes C_4$. For the center of this group we get 
$Z(\Aut(M_2)) \cong C_2$ and its system of normal subgroups can be seen in 
the next figure.
\begin{center}
\setlength{\unitlength}{1cm}
\begin{picture}(5,5)
\put(2.4,4.7){id}
\put(1.5,3.5){\line(1,1){1.0}}
\put(3.5,3.5){\line(-1,1){1.0}}
\put(1.3,3.0){$C_2$}
\put(3.4,3.0){$A_5$}
\put(2.5,1.8){\line(-1,1){1.0}}
\put(2.5,1.8){\line(1,1){1.0}}
\put(1.88,1.4){$C_2 \times A_5$}
\put(2.5,0.2){\line(0,1){1.0}}
\put(1.85,-0.2){$A_5 \rtimes C_4$}
\vspace{2mm}
\end{picture} 
\end{center}
\begin{center}{\bf Figure 2.} Normal subgroups of $\Aut(M_2)$.
\vspace{2mm}
\end{center}

Furthermore it holds
$\Aut(M_2) / Z(\Aut(M_2)) \cong S_5$, but $S_5$ is not a subgroup of
$\Aut(M_2)$.
\end{proof}
\vspace{2mm}
Note that Lemma $2.3$ immediately yields 
$M_1 \not\sim M_2 \not\sim M_3 \not\sim M_1$.\\[2mm]

\begin{Rem}
The program {\em \texttt{SUNI}} has found the following three other
candidates:
\renewcommand{\labelenumi}{(\arabic{enumi})} 
\begin{enumerate}
\setcounter{enumi}{3}
\item
$ \quad 1\, 1\, 1\, 1\,8 \qquad 2\,2\,3\,2\,3 \qquad 1\,2\,2\,3\,4 \qquad
1\,1\,7\,1\,2 \qquad 1\,1\,2\,5\,3 \qquad 1\,3\,1\,3\,4$
\item 
$ \quad 1\, 1\, 1\, 1\,8 \qquad 2\,2\,3\,2\,3 \qquad 1\,2\,2\,3\,4 \qquad
1\,1\,7\,1\,2 \qquad 1\,1\,2\,5\,3 \qquad 1\,3\,1\,4\,3$
\item
$ \quad 1\, 1\, 2\, 1\,7 \qquad 1\,4\,3\,2\,2 \qquad 1\,2\,2\,3\,4 \qquad
1\,1\,7\,1\,2 \qquad 1\,1\,2\,5\,3 \qquad 1\,1\,3\,5\,2$
\vspace{2mm}
\end{enumerate}
\renewcommand{\labelenumi}{(\roman{enumi})} 
We have $M_1 \sim M_5$ and $M_2 \sim M_4$ as can be seen immediately be
applying the multiplier $-1 \equiv 11 \mod 12$ on $M_1$, $M_2$ respectively.
(For the notion of multipliers see Section $3$.)\\
This can also be proved by the {\em \texttt{GAP}}-program mentioned in the 
proof of $2.2$. Any of the calculated bijections between the vertex links of 
two examples can be extended to give a bijection between the examples 
themselves.\\

Example $M_6$ is not a manifold as we prove now:\\

As $\langle 048 \rangle \not\in M_6$, $M_6$ is not $2$-Hamiltonian in 
$C_6^\ast$ and therefore not a triangulated manifold.
(An upgrade of the program {\em \texttt{SUNI}} meanwhile also calculates the 
$f$-vector of all candidates. The output of this extension does not 
contain $M_6$ anymore.)
\vspace{6mm}

Alternatively, one can also examine the link 
$lk(\langle 04 \rangle)$ of 
edge $\langle 04 \rangle$. It consists of the following $2$-simplices:
\[\begin{array}{cccccc}  \langle 125\rangle, & 
\langle 13b\rangle, & \langle 79b\rangle, & 
\langle 259\rangle, & \langle 27b\rangle, & 
\langle 579\rangle, \\ \langle 135\rangle, & 
\langle 23b\rangle, & \langle 129\rangle, & 
\langle 357\rangle, & \langle 19b\rangle, & 
\langle 237\rangle.
\end{array} \]
Therefore $f_0(lk(\langle 04 \rangle)) = 7$, $f_1(lk(\langle 04 \rangle)) = 
18$, $f_2(lk(\langle 04 \rangle)) = 12$ and $\chi(lk(\langle 04 \rangle)) = 
1$. Hence $M_6$ is not even an Eulerian $4$-manifold, i.~e.\ 
$\chi(lk(\Delta^k)) = 1 -(-1)^{4-k}=1-(-1)^k$ is not satisfied for all 
$k$-simplices of $M_6$.

But $lk(\langle 04 \rangle)$ is a pinched $2$-sphere and therefore a
$2$-manifold with one singularity as can be seen from Figure $3$.
\begin{center}
\epsfig{file=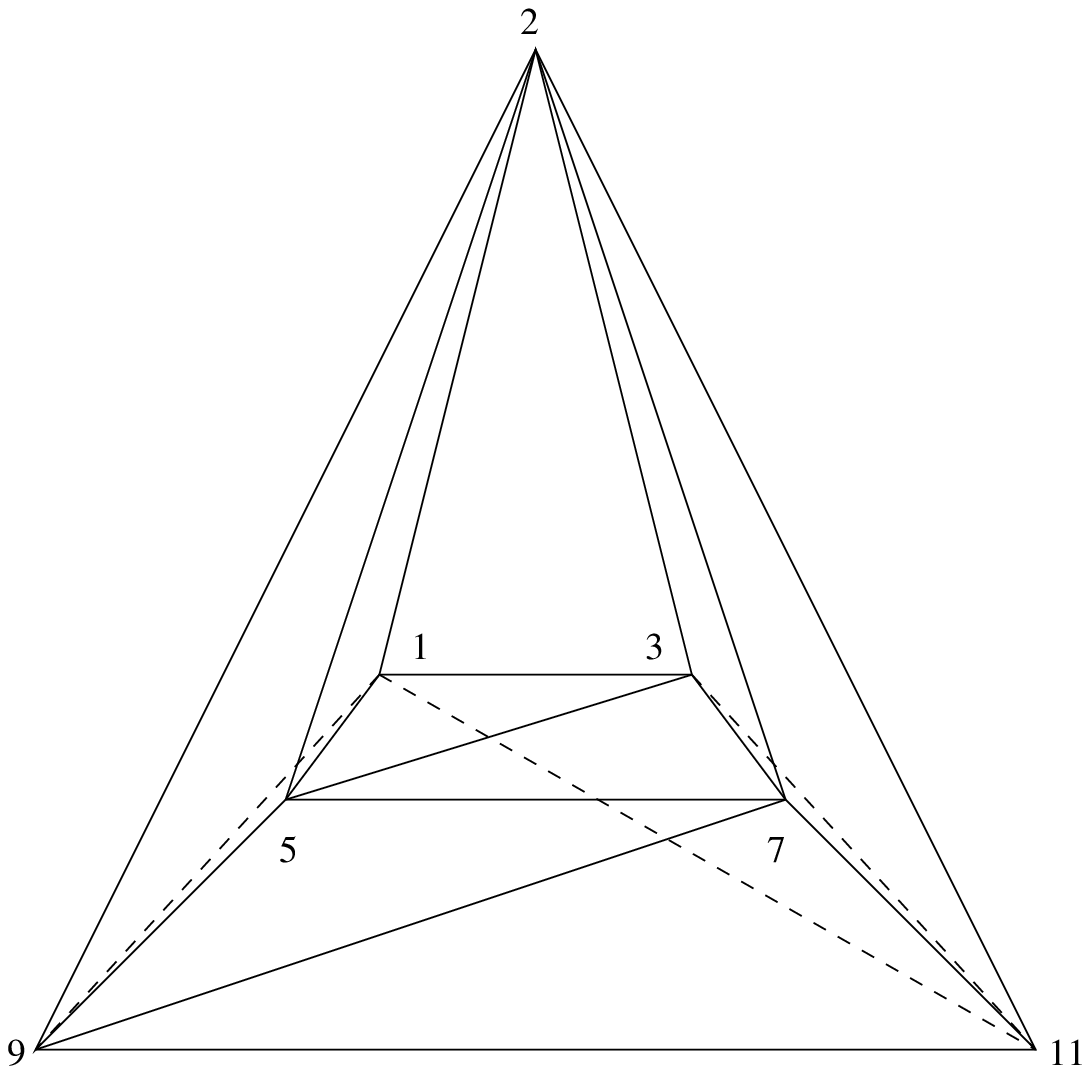, width=9.4cm, height=5.2cm} \\
{\bf Figure 3.} $lk(\langle 04 \rangle)$ of example $M_6$
\vspace{5mm} 
\end{center}
\vspace{-3mm}
\end{Rem}
\vspace{3mm}

\section{Properties of the triangulations of $S^2\times S^2$}
In this section we describe further properties of the examples $M_i$. 
We begin with a simple observation.\\

The span of the ``even'' vertices as well as the span of the ``odd'' 
vertices of all $M_i$ 
is a subcomplex of the octahedron, i.~e.\ the $3$-cross-polytope.
This holds because e.~g.\ $span(0,2,4,6,8,a)$ contains no $k$-simplex with $k
\geqslant 3$. (Otherwise there would be a $3$-simplex in $span(0,2,4,6,8,a)$ 
containing a diagonal as an edge.) The $f$-vector of the $M_i$ implies, that in
 $span(0,2,4,6,8,a)$ all $2$-simplices without one of the diagonals
$\langle 06\rangle$, $\langle 28\rangle$ and $\langle 4a\rangle$ are
contained. Therefore we get that  $span(0,2,4,6,8,a)$ equals the boundary 
complex of $C_3^\ast$. Analogously we get the same result for the odd vertices.

This not only holds for the even and odd vertices, but for all
subsets $S$ of the vertex set with $|S|=6$, such that
$3$ diagonals are contained in the regarded subset.\\

Another interesting subset is the cylinder $C$ consisting of the triangles
\[ \langle 024\rangle,\;  \langle 246\rangle,\;  \langle 468\rangle,\; 
   \langle 68a\rangle,\;  \langle 8a0\rangle \text{ and }  
   \langle 0a2\rangle. 
\]
This cylinder is invariant under automorphism $\zeta^2$. Together with
the triangles $\langle 048\rangle$ and $\langle 26a\rangle$, whose
boundaries form the boundary of $C$, we get $span(0,2,4,6,8,a)$. \\
Likewise there is a cylinder in $span(1,3,5,7,9,b)$.\\[2mm]

Now let us regard $lk(0)$ of example $M_2$. The vertices of this link  
can be subdivided into two subsets $S_{2_1}$ and $S_{2_2}$ 
with each $5$ elements, such that the span of both subsets form M\"obius 
strips in $lk(0)$ (cf.\ \cite{Sp1}). As all vertex links are combinatorally 
equivalent, the same holds for any vertex link.

It is easy to show that the required subdivision is unique:
Because any $5$-vertex M\"obius strip must be neighborly, no two elements of 
$S_{2_1}$ may form a diagonal of $M_2$. Suppose without loss of generality 
$1 \in S_{2_1}$. Then $7 \in S_{2_2}$ and there remain $8\cdot 6 \cdot 4 
\cdot 2 = 384$ possibilities to check.\\
In \cite{Sp1} we used these M\"obius strips to compute the 
intersection numbers. In this case however it was not necessary that we got
M\"obius strips, but that the M\"obius strips were linked and unknotted.
So although there is no subdivision of any vertex link of either $M_1$ or 
$M_3$, 
that also produces disjoint M\"obius strips, the proof still carries over to 
the cases $M_1$ and $M_3$. Nonetheless the existence of disjoint M\"obius
strips stills ``shows'' that example $M_2$ is the most symmetric one.\\[2mm]

Inequality $(1)$ gives an upper bound on $|\Aut(M_i)|$. We will derive now a
lower bound. For this we use the notion of multipliers. 
\begin{Def}
Let the vertices of a simplical complex $K$ be numbered by $0,\ldots, 
\linebreak k-1$.
$a \in C_k$ is called a multiplier of $K$ if 
\[ \mu_k : C_k \rightarrow C_k,\quad i \mapsto i \cdot a \mod k\]
is an automorphism of $K$.\\
\end{Def}

Recall the following simple facts:
\begin{enumerate}
\item If $a$ is a multiplier of $K$, then $gcd(a,k)=1$.
\item If $(0,\ldots,k-1)$ is an automorphism of $K$ and if $K$ has $l$ 
different multipliers, then $|\Aut(K)| \geqslant l\cdot k$.\\
\end{enumerate}

From $(i)$ we know that the only possible multipliers of our 
triangulations are $1,5,7$ and $11$. $M_1$ and $M_2$ have multipliers $1$ 
and $5$, whereas under the maps $\mu_7$ and $\mu_{11}$ the image of 
$M_1$ is $M_2$ and vice versa. $M_3$ has all multipliers. Hence by $(ii)$
\[ |\Aut(M_i)|  \geqslant \begin{cases} 
24 & \text{for } i=1,\;2\\ 48 & \text{for } i=3\end{cases}\;.
\]
Together with $(1)$ we conclude  
\[ |\Aut(M_i)|  = \begin{cases} 
24 & \text{for } i=1\\ 48 & \text{for } i=3\end{cases}\;.
\]
\vspace{2mm} \relax{} \\\

\section{Further $12$- and $11$-vertex triangulations}

In the next remark we describe another interesting simplical complex
found by \texttt{SUNI}.

\begin{Rem} The program {\em \texttt{SUNI}} has delivered another candidate $M$ for being a 
\linebreak $4$-manifold with $12$ vertices and with cyclic symmetry. In this 
example we dropped the 
additional assumption that $\langle 06\rangle$ has to be a diagonal. 
This example is rather interesting as it is highly-symmetric, as one can see
from its invariance under {\em all} multipliers. It is given by the 
difference-$5$-cycles
\[ 1\,1\,1\,3\,6 \quad 1\,4\,2\,3\,2 \quad 1\,2\,3\,2\,4 \quad 
1\,1\,1\,6\,3 \quad 1\,1\,3\,3\,4 \quad 1\,3\,1\,5\,2 \quad 
1\,2\,5\,1\,3 \quad 1\,1\,4\,3\,3 \]
and we get for the $f$-vector $(f_1,\ldots,f_4)=(12,66,204,240,96)$ and
consequently $\chi(M)=6$. (Missing triangles are for instance 
$\langle 024\rangle$ and $\langle 048\rangle$; cp.\ the properties of 
$M_1$ - $M_3$ as described on the previous page.)
However this candidate is not a \linebreak $4$-manifold, as the link of 
the edge 
$\langle 01\rangle$ consists of three $2$-spheres that are only edge- but not
vertex-disjoint. These spheres are clued together at their poles
(``Banana-surface'') as shown in the next figure.  
\begin{center}
\epsfig{file=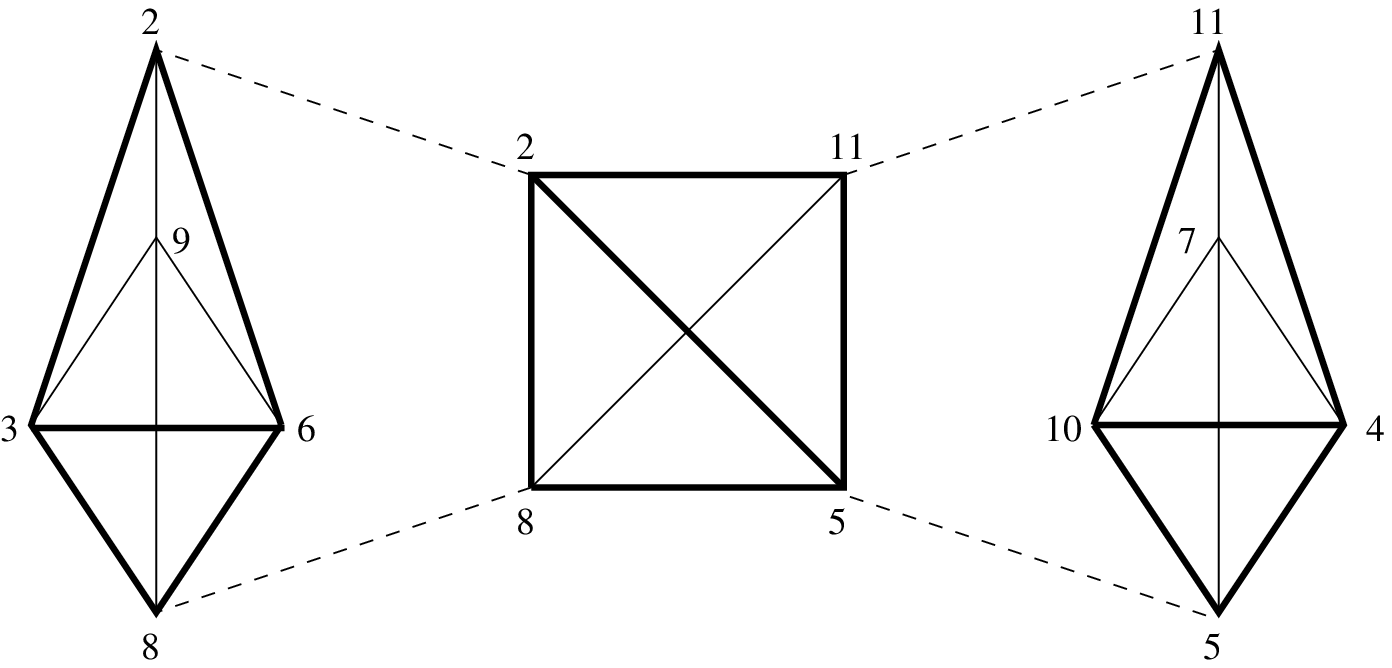, width=10cm} \\
{\bf Figure 3.} The link of edge $\langle 0 1 \rangle$. \vspace{3mm}
\end{center}
But $M$ is at least an Eulerian $4$-manifold  \cite{Ku2} as follows from 
calculating the Euler-characteristic of the links of all simplices $\Delta^k$
of the triangulation. Because of the cyclic symmetry it suffices to show
$\chi(lk(\langle i \rangle)) = 0$ for one vertex. {\em \texttt{SUNI}} 
\linebreak itself
confirms $\chi(lk(\Delta^1)) = 2$ for all $1$-simplices of the triangulation
and finally $\chi(lk(\Delta^2)) = 0 = 2-\chi(lk(\Delta^3))$ can easily be 
checked.\\
The singular locus of this example is the union of two tori. One of these
is given by the difference-$3$-cycles $012$ and $05a$ $\!\!\mod\!$ $12$, the other by
$015$ and $045$. Each of these tori is invariant under the group of order
$48$ generated by the cyclic automorphism
$\zeta$ and the four multipliers $\mu_k$.\\
\end{Rem}

\begin{Rem}
Homology investigations of Frank Lutz (oral communication) imply that their 
are no other $12$-vertex triangulations of $S^2 \times S^2$ with an 
automorphism group that is transitive on the vertices than the examples 
found in our Main \linebreak Theorem. As these are all rather symmetric (especially 
example $M_2$), we felt that it will be possible
to ``break'' a bit of the symmetry of at least one of our examples to
construct a triangulation of $S^2 \times S^2$ with just $11$ vertices.

That the theoretical lower bound on the number of vertices according to  
\cite{Ku1} is $10$ and that this value was ruled out only by an extensive 
computer search made this assumption even more plausible.

While preparing this paper an unpublished {\em \texttt{GAP}}-program written 
by Frank Lutz shows that there really is such an example:\\

{\em There does exist a triangulation of $S^2 \times S^2$ with $11$ 
vertices -- the minimum number of vertices.}\\

It remains open if this $11$-vertex-triangulation is unique. We conjecture 
that this is not the case. Using bistellar operations one can possibly
construct $S^2 \times S^2$ with $11$ vertices not combinatorially
equivalent to the given triangulation.\\

Apart from a computer search as applied by Frank Lutz, a more systematic way 
of constructing such a triangulation may be to start with one 
of our examples $M_i$ and cut off the open star of one vertex. The resulting 
complex $\widetilde{M}_i$ has $11$ vertices and is homeomorphic to 
$S^2 \times S^2$ with an open ball 
removed. Now try to close $\widetilde{M}_i$ by inserting new simplices, but no
new vertices. Although this procedure may not always be possible, it can often 
be carried out \cite{ABS}. Here we may use simplices that 
contain edges that are diagonals in $M_i$; this is the process we refered 
to as ``breaking of 
symmetry''. Finding such simplices naturally leads a system of linear
equations.
\end{Rem}

\end{document}